\documentclass[12pt]{amsart}

\usepackage{amssymb, amsthm, amsfonts, amsmath}

\newcommand{\A}{\mathcal{A}}
\newcommand{\B}{\mathcal{B}}
\newcommand{\C}{\mathcal{C}}
\newcommand{\DD}{\mathcal{D}}
\newcommand{\EE}{\mathcal{E}}
\newcommand{\FF}{\mathcal{F}}
\newcommand{\GG}{\mathcal{G}}
\newcommand{\KK}{\mathcal{K}}
\newcommand{\LL}{\mathcal{L}}
\newcommand{\MM}{\mathcal{M}}
\newcommand{\OO}{\mathcal{O}}
\newcommand{\into}{\hookrightarrow}

\newcommand{\mm}{\mathfrak{m}}
\def\from{\leftarrow}
\def\spot{{\hbox{\raisebox{.33ex}{\large\bf .}}\hspace{-.10ex}}}

\DeclareMathOperator{\Spec}{Spec}
\DeclareMathOperator{\TorTor}{\ensuremath{\mathcal{T}\hspace{-.7ex}\mathit{or}}}
\DeclareMathOperator{\Tor}{Tor}

\newtheorem*{MT}{Theorem}
\newtheorem*{GFT}{Generic Freeness}
\newtheorem{prop}{Proposition}
\newtheorem{lemma}[prop]{Lemma}

\theoremstyle{definition}
\newtheorem*{remark}{Remark}

\newenvironment{numbered}%
        {\begin{list}
                {\noindent\makebox[0mm][r]{\arabic{enumi}.}}
                {\leftmargin=5.5ex \usecounter{enumi}}
        }
        {\end{list}}

\begin{document}

\mbox{}
\vspace{-5.6ex}
\title{A Kleiman--Bertini Theorem for sheaf tensor products}
\author{Ezra Miller}
\address{School of Mathematics, University of Minnesota, Minneapolis, MN, USA}
\email{ezra@math.umn.edu\vspace{-1.5ex}}

\author{David E Speyer}
\address{Department of Mathematics, University of Michigan, Ann Arbor, MI, USA}
\email{speyer@umich.edu}

\thanks{EM gratefully acknowledges support from NSF CAREER award
DMS-0449102 and a University of Minnesota McKnight Land-Grant
Professorship.  DES is a Clay Research Fellow, and is pleased to
acknowledge the support of the Clay Mathematics Institute.  This paper
originated in a visit of DES to the University of Minnesota, and he is
grateful for their excellent hospitality.}

\date{8 February 2007}

\begin{abstract}
Fix a variety $X$ with a transitive (left) action by an algebraic
group~$G$.  Let $\EE$ and $\FF$ be coherent sheaves on~$X$.  We prove
that for elements $g$ in a dense open subset of~$G$, the sheaf
$\TorTor^X_i(\EE, g \FF)$ vanishes for all $i > 0$.  When $\EE$ and
$\FF$ are structure sheaves of smooth subschemes of $X$ in
characteristic zero, this follows from the Kleiman--Bertini theorem;
our result has no smoothness hypotheses on the supports of $\EE$
or~$\FF$, or hypotheses on the characteristic of the ground field.
\end{abstract}
  
\maketitle

All schemes in this note are of finite type over an arbitrary base
field~$k$.  We make no restrictions on the characteristic of~$k$.  All
groups~$G$ are assumed to be smooth over~$k$.  In particular, the
local rings of $G \times_k k'$ are regular, for all field
extensions~$k'$ of~$k$.  By a \emph{transitive} action of~$G$ on~$X$,
we mean one such that $G \times X \to X \times X$ is
scheme-theoretically surjective (this implies that $X$ is reduced; see
the Remark).

The aim of this note is to prove the following result.

\begin{MT}
  Let $X$ be a variety with a transitive left action of an algebraic
  group~$G$.  Let $\EE$ and $\FF$ be coherent sheaves on $X$ and, for
  all $k$-rational points $g \in G$, let $g \FF$ denote the
  pushforward of $\FF$ along multiplication by $g$.  Then there is a
  dense Zariski open subset~$U$ of~$G$ such that, for $g \in U$, the
  sheaf $\TorTor_i(\EE, g \FF)$ is $0$ for all~\mbox{$i > 0$}.
\end{MT}

\begin{remark}
There are a number of subtleties that can be ignored for those working
over the field of complex numbers.
\begin{numbered}
\item
Transitivity implies that the action is \emph{geometrically
transitive}, in the sense that, after base extending to an algebraic
closure of~$k$, any choice of closed point in~$X$ results in a
scheme-theoretically surjective map $G \to X$.  Hence $X$ is
geometrically reduced, and therefore reduced.

\item
If $k$ is infinite, and $G$ is connected as well as affine, and either
(i)~$G$ is reductive or (ii)~$k$ is perfect, then $U$ always contains
a $k$-rational point of~$G$ \cite[Corollary~18.3]{Bor}.

\item
For the open set $U$ that we construct, if $k'$ is any extension
of~$k$ and $g$ is a $k'$-rational point of $U$, then
$\TorTor_i^{X'}(\EE, g \FF)$ will vanish as well, where $X' = X
\times_k k'$ is the base extension to~$k'$.
\end{numbered}
\end{remark}

The theorem should be thought of as an analogue of the
Kleiman--Bertini theorem~\cite{Klein}, which concerns the
characteristic zero case in which $\EE = \OO_Y$ and $\FF = \OO_Z$ are
structure sheaves of smooth subvarieties: $Y$ meets $g Z$
transversally for generic $g \in G$ in that situation, immediately
implying that $\EE$ and $g \FF$ have no higher $\TorTor$.  The
conclusion that $Y$ and $Z$ meet transversally has many additional
consequences; it implies, for example, that $Y \cap g Z$ is smooth,
which is the better known conclusion of the Kleiman--Bertini theorem.
In our situation, where $\EE$ and $\FF$ need not be structure sheaves,
there is no analogous smoothness conclusion to draw.  Our result shows
that we may at least continue to disregard higher $\TorTor$ sheaves in
this general setting.

There are many results which generalize the Kleiman--Bertini theorem
by saying that for $X$ as in our theorem, if $Y$ and $Z$ are
subvarieties of $X$ obeying various conditions on their singularities,
then for generic $g \in G$ the intersection $Y \cap g Z$ obeys the
same conditions (sometimes we must assume that the ground field has
characteristic zero).  The cases where the condition on singularities
is reducedness (in characteristic zero), equidimensionality, local
connectivity in codimension $c$, Serre's conditions $R_d$ (in
characteristic zero) or $S_d$, or normality can be obtained by methods
similar to those of Spreafico \cite{Sprea}.  (Spreafico considers
Kleiman--Bertini-like results in Section~5.2 of her paper.  In order
to avoid imposing Spreafico's condition~(b), which is not true for
many interesting examples such as grassmannians, one should remove her
notion of generic geometric regularity and instead impose on~$P$ the
axioms that (1)~if $A$ and~$B$ have property~$P$ then $A \times B$ has
property~$P$, and (2)~if $A$ has property~$P$, and $B$ is a regular
domain with a morphism $A \to B$, then $A \times_B \overline\eta$ has
property~$P$, where $\overline\eta$ is a geometric generic point
of~$B$.  It is checking this second property that introduces
characteristic dependencies.)  The condition of rationality of
singularities is essentially Lemmas~1 and~2 of \cite{Brion}.  These
lemmas also provide an alternative proof when the condition on
singularities is Cohen--Macaulayness, and it is straightforward to
modify this argument to apply to Gorenstein singularities.  In the
same paper, Brion also proves the special case of our result where
$\EE$ and~$\FF$ are the structure sheaves of Cohen--Macaulay
subschemes of~$X$.  We could not find in the literature a result
dealing with sheaves that are not structure sheaves.

The main use we see for our result is in $K$-theory computations for
homogeneous spaces.  This is an active field of research, with many
results both in pure geometry and in combinatorics; see \cite{Brion2}
for a good introduction to the geometric side of the theory and
\cite{Buch02,Las90} for a sample of the combinatorial side.  The
$K$-ring $K^\circ(X)$ of a smooth variety~$X$ can be defined
additively to be the abelian group generated by the collection of
coherent sheaves on~$X$ subject to the relations $[\A]+[\C]=[\B]$
whenever there is a short exact sequence $0 \to \A \to \B \to C \to
0$.  The multiplication in $K^\circ(X)$ is induced by the tensor
product; more precisely, it is given by
\begin{equation} \label{KProduct}
[\EE][\FF]=\sum_i (-1)^i [\TorTor^X_i(\EE,\FF)].
\end{equation}
(The smoothness of~$X$ guarantees that this sum terminates.)  Whenever
$G$ is a connected linear algebraic group acting transitively on $X$,
the $G$-action on~$K^\circ(X)$ is trivial, so we may compute $[\EE][g
\FF]$ instead of~$[\EE][\FF]$.  Our result says that when $g$ is
chosen generically, only the $i = 0$ term in~(\ref{KProduct}) can be
nonzero.  This is useful, as the tensor product is a far more familiar
object than the higher $\TorTor$ sheaves.
In particular, when $\EE = \OO_Y$ and $\FF = \OO_Z$ are structure
sheaves for arbitrary closed subschemes of~$X$, the tensor product
$\EE \otimes \FF$ is simply the structure sheaf $\OO_{Y \cap Z}$ of
the scheme-theoretic intersection $Y \cap Z$.

\section*{Proof of Theorem}

{}From now on, $G$, $X$, $\EE$ and $\FF$ have the meanings in the
statement of our Theorem.  Let $\Gamma$ denote the subvariety of $X
\times X \times G$ consisting of those points $(x_1,x_2,g)$ such that
$x_1 = g x_2$.  Thus $\Gamma$ is an isomorphic copy of $X \times G$
embedded in $X \times (X \times G)$ by the graph of the multiplication
map $X \times G \to X$.  We will write $p_1$, $p_2$, and $q$ for the
projections of~$\Gamma$ onto the two $X$ factors and~$G$,
respectively.

\begin{lemma} \label{flat}
The map $p_1 \times p_2 : \Gamma \to X \times X$ is flat.
\end{lemma}
\begin{proof}
Extending the base field~$k$, we assume $X$ has a $k$-rational point,
with stabilizer~$H$.  (Note that $H$ need not be reduced; this causes
no difficulty here or elsewhere.)  By virtue of the transitivity of
the group action, $X$ is smooth, the map $\mu:G \to X$ is surjective,
and all fibers of $\mu$ have the same dimension. By our standing
assumption, $G$ is smooth, so the map $\mu$ is faithfully flat by
\cite[Proposition~15.4.2]{EGA4}.  Base extending $p_1 \times p_2$ by
$G \times G \to X \times X$, we need only prove that $G \times G
\times_{X \times X} \Gamma \to G \times G$ is flat.  Thus, if $\Delta
= \{(g_1, g_2, g) : g g_2 \in g_1 H\} \subseteq G \times G \times G$,
we need the projection $p_{12}$ of~$\Delta$ to the first two factors
to be flat.  But the automorphism of $G \times G \times G$ taking
$(g_1, g_2, g) \mapsto (g_1, g_2, g_1^{-1} g g_2)$ commutes
with~$p_{12}$ and takes~$\Delta$ to $G \times G \times H$.%
\end{proof}

Lemma~\ref{flat} would be false without the reduced (no nilpotents)
hypothesis on~$X$, as can be seen by taking $X =
\Spec(k[\varepsilon]/\varepsilon^2)$ and setting $G = \Spec(k)$, the
trivial group.

Let $\GG$ be the coherent sheaf $p_1^*\EE \otimes p_2^*\FF$
on~$\Gamma$.  Thus $\GG = (p_1 \times p_2)^*(\EE \boxtimes \FF)$,
where $\boxtimes$ denotes the tensor product of the pullbacks to $X
\times X$ along its two projections to~$X$.  For $g \in G(k)$, let
$\iota^g : X \into \Gamma$ be the map $x \mapsto (x, g^{-1}x, g)$.
Let $k_g$ be the skyscraper sheaf on $G$ concentrated at $g$.  The
heart of our proof is the following computation.

\begin{prop} \label{differentTors}
With the above notation, we have
\[
  \iota^g_* \TorTor^X_i(\EE, g \FF) = \TorTor^\Gamma_i(q^* k_g, \GG).
\]
\end{prop}
\begin{proof}
$X$ is smooth, by virtue of the transitive group action, so $X$ has
enough locally frees (see Exercises~6.4, 6.8, and~6.9 in
\cite[Chapter~III]{Har}).  Let $0 \from \KK_0 \from \KK_1 \from
\cdots$ and $0 \from \LL_0 \from \LL_1 \from \cdots$ be resolutions of
$\EE$ and~$\FF$ by locally free $\OO_X$-modules.  Then $\C_{i,j}=\KK_i
\boxtimes \LL_j$, is a double complex on $X \times X$.  The rows
$\C_{\spot,j}$ resolve $\EE \boxtimes \LL_j$, the columns
$\C_{i,\spot}$ resolve $\KK_i \boxtimes \FF$, and the total complex
resolves $\EE \boxtimes \FF$.  Let
\[
  \DD_{i,j} = (p_1 \times p_2)^* (\KK_i \boxtimes \LL_j) = p_1^*\KK_i
  \otimes p_2^*\LL_j.
\]
Because $p_1 \times p_2$ is a flat morphism, the rows, columns, and
total complex of $\DD_{\spot,\spot}$ resolve \mbox{$p_1^*\EE \otimes
p_2^*\LL_j$}, $p_1^*\KK_i \otimes p_2^*\FF$, and $\GG = p_1^*\EE
\otimes p_2^*\FF$, respectively.

We claim that the rows of $\DD_{\spot,\spot} \otimes q^*k_g$ can only
have nonzero homology in homological degree~$0$, and that in row~$j$
this homology is $p_1^*\EE \otimes p_2^*\LL_j \otimes q^*k_g$.  To see
this, note that $p_1^*\MM \otimes q^*k_g = \iota^g_*\MM$ for all
coherent sheaves~$\MM$ on~$X$.  Since $\iota^g$ is a closed embedding,
$\iota^g_*$ is exact, so the complex $0 \from p_1^*\KK_0 \otimes
q^*k_g \from p_1^*\KK_1 \otimes q^*k_g \from \cdots$ is a resolution
of $\iota^g_*\EE = p_1^*\EE \otimes q^*k_g$.  Tensoring this
resolution with the locally free sheaf~$p_2^*\LL_j$ preserves
exactness, resulting in a resolution of $p_1^*\EE \otimes p_2^*\LL_j
\otimes q^*k_g$, as~claimed.
  
Now, consider the horizontal homology of $\DD_{\spot,\spot} \otimes
q^*k_g$, which is concentrated in the zeroth column, as a complex
under the vertical differential, namely
\begin{equation} \label{blah}
  0\ \from\ p_1^*\EE \otimes p_2^*\LL_0 \otimes q^*k_g\ \from\
  p_1^*\EE \otimes p_2^*\LL_1 \otimes q^*k_g\ \from\ \cdots.
\end{equation}
By an argument similar to that of the previous paragraph, we may
rewrite~(\ref{blah}) as the pushforward of a complex on~$X$; namely,
(\ref{blah}) is the image under $\iota^g_*$ of
\begin{equation} \label{blah2}
  0\ \from\ \EE \otimes g\LL_0\ \from\ \EE \otimes g\LL_1\ \from\
  \cdots,
\end{equation}
where we have used that $p_2^*\MM \otimes q^*k_g = \iota^g_*(g\MM)$
for all coherent sheaves~$\MM$ on~$X$.  The $i^{\textrm{th}}$ homology
of~(\ref{blah2}) is $\TorTor^X_i(\EE, g \FF)$.  As $\iota^g$ is a
closed immersion, $\iota^g_*$ is exact; hence the $i^{\textrm{th}}$
homology of~(\ref{blah}) is $\iota^g_* \TorTor^X_i(\EE, g \FF)$.

To summarize, $\DD_{\spot,\spot} \otimes q^*k_g$ is a double complex
whose rows are acyclic and whose horizontal homology is a vertical
complex computing $\iota^g_* \TorTor^X_i(\EE, g \FF)$.  On the other
hand, the total complex of $\DD_{\spot,\spot}$ is a locally free
resolution of~$\GG$ over~$\OO_\Gamma$, so the $i^{\textrm{th}}$
homology of~$\DD_{\spot,\spot} \otimes q^*k_g$ is
$\TorTor^\Gamma_i(q^*k_g, \GG)$.  Standard homological techniques (see
\cite[Section~5.6]{Hom}, for example) produce a natural isomorphism
from the homology of the total complex to that of the vertical
complex, given the horizontal acyclicity.%
\end{proof}

Due to the isomorphism in Proposition~\ref{differentTors}, we are
interested in $\TorTor^\Gamma_i(q^* k_g, \GG)$.
\begin{prop} \label{genericVanishing}
For $i>0$ and generic $g \in G$, we have $\TorTor^\Gamma_i(q^*k_g,\GG)
= 0$.
\end{prop}

The proof is based on the following general result due to
Grothendieck.

\begin{GFT}
Let $A$ be a generically reduced noetherian scheme, $q: B \to A$ a
finite type $A$-scheme, and $\MM$ a coherent sheaf on~$B$.  Then there
is a dense open subset~$U$ of~$A$ such that (the pushforward to~$U$
of) $\MM|_{q^{-1}(U)}$ is a locally free~\mbox{$\OO_U$-module}.
\end{GFT}
\begin{proof}
See \cite[$\textrm{IV}_3$,~11.2.6.1]{EGA3}, or see
\cite[Theorem~14.4]{Eis} for a very readable proof in the affine case.
\end{proof}

\begin{proof}[Proof of Proposition~\ref{genericVanishing}]
Apply generic freeness twice to find a dense open $U \subset G$ such
that $\GG|_{q^{-1}(U)}$ and $\OO_{\Gamma}|_{q^{-1}(U)}$ are both
locally free $\OO_U$ modules.  We claim that
$\TorTor^\Gamma_i(q^*k_g,\GG) = 0$ for all $i>0$ when $g \in U$.  This
computation may be checked locally: if $R$ is the local ring of
$\OO_G$ at $g \in U$, with maximal ideal~$\mm_g$, and $S$ is the local
ring of any point $\gamma \in \Gamma$ mapping to~$g$, then we only
need the $S$-module $\Tor^S_i(S \otimes_R R/\mm_g,M)$ to vanish for $i
> 0$, where $M = \GG_\gamma$ is the stalk of~$\GG$ at~$\gamma$.

Considering $\Tor^S_i(S \otimes_R R/\mm_g,M)$ as an $R$-module, it
agrees with $\Tor^R_i(R/\mm_g,M)$, since any resolution~$F_\spot$
of~$M$ by free $S$-modules is already a resolution of~$M$ by free
$R$-modules (by our choice of $U$), and $F_\spot \otimes_S (S \otimes_R
R/\mm_g) = F_\spot \otimes_R (R/\mm_g)$.  The desired vanishing
follows because $M$ is a free (hence flat) $R$-module (by our choice of $U$).%
\end{proof}

The proof of our Theorem is now easy to complete.
Propositions~\ref{differentTors} and~\ref{genericVanishing} imply that
for $g \in G$ generic, $\iota^g_*\TorTor^X_i(\EE, g\FF) =
\TorTor^\Gamma_i(q^*k_g,\GG) = 0$ whenever $i > 0$.  But $\iota^g$ is
a closed immersion, so $\iota^g_*$ is exact and faithful; in
particular, $\iota^g_*\TorTor^X_i(\EE, g\FF) = 0$ implies that
$\TorTor^X_i(\EE, g\FF) = 0$, as desired.  \nolinebreak\qedsymbol

\bigskip
\noindent
\textbf{Acknowledgements.}  The authors are grateful to Brian Conrad,
Igor Dolgachev, and Michel Raynaud for helpful suggestions and
references.

\vspace{-.7ex}

\end{document}